\documentclass[a4paper,10pt,leqno,twoside]{article}

\usepackage[english]{babel} \usepackage{inputenc, amsmath, amssymb , latexsym,
  epic, epsfig, rotating, fancyheadings, amsthm, pifont, empheq}

\newcommand{\homd}{\ensuremath{\mathrm{Homeo}_0(\mathbb{T}^d)}}

\newcommand{\follows}{\ensuremath{\Rightarrow}}

\newcommand{\ld}{\ensuremath{,\ldots,}}
\newcommand{\ssq}{\ensuremath{\subseteq}}
\newcommand{\smin}{\ensuremath{\setminus}}
\newcommand{\eps}{\ensuremath{\varepsilon}}

\newcommand{\T}{\ensuremath{\mathbb{T}}}

\newcommand{\N}{\ensuremath{\mathbb{N}}} 
\newcommand{\R}{\ensuremath{\mathbb{R}}}
\newcommand{\Z}{\ensuremath{\mathbb{Z}}}
\newcommand{\Q}{\ensuremath{\mathbb{Q}}}

\newcommand{\sph}{\ensuremath{\mathbb{S}}}

\newcommand{\Id}{\ensuremath{\mathrm{Id}}}

\newcommand{\supp}{\ensuremath{\mathrm{supp}}}

\newcommand{\Leb}{\ensuremath{\mathrm{Leb}}}

\newcommand{\inte}{\ensuremath{\mathrm{int}}}

\newcommand{\kreis}{\ensuremath{\mathbb{T}^{1}}}
\newcommand{\ntorus}[1][2]{\ensuremath{\mathbb{T}^{#1}}}

\newcommand{\torus}{\ensuremath{\mathbb{T}^2}}

\newcommand{\dhom}{\ensuremath{d_\mathrm{hom}}}

\newcommand{\alphlist}{\begin{list}{(\alph{enumi})}{\usecounter{enumi}}}
\newcommand{\romanlist}{\begin{list}{(\roman{enumi})}{\usecounter{enumi}}}
\newcommand{\listend}{\end{list}}

\newcommand{\roundqed}{{\raggedleft \Large $\circ$\\}}

\newcommand{\nfolge}[1]{\ensuremath{(#1)_{n\in\mathbb{N}}}}

\newcommand{\incup}{\ensuremath{\bigcup_{i=1}^n}}

\newcommand{\kcap}{\ensuremath{\bigcap_{k\in\N}}}

\newcommand{\nLim}{\ensuremath{\lim_{n\rightarrow\infty}}}
\newcommand{\iLim}{\ensuremath{\lim_{i\rightarrow\infty}}}
\newcommand{\jLim}{\ensuremath{\lim_{j\rightarrow\infty}}}

\newcommand{\ntel}{\ensuremath{\frac{1}{n}}}

\newcommand{\halb}{\ensuremath{\frac{1}{2}}}

\newcommand{\achtel}{\ensuremath{\frac{1}{8}}}

\newtheorem{definition}{Definition}[section]
\newtheorem{thm}[definition]{Theorem}

\newtheorem{lem}[definition]{Lemma}
\newtheorem{cor}[definition]{Corollary}
\newtheorem{prop}[definition]{Proposition}
\newtheorem{claim}[definition]{Claim}
\newtheorem{question}[definition]{Question}
\newtheorem{bem}[definition]{Remark}
\newtheorem{example}[definition]{Example}

\numberwithin{equation}{section}

\title{\Large\textsc{The concept of bounded mean motion for toral
    homeomorphisms}} \author{T.~J\"ager\thanks{Coll\`ege de France,
    Paris. Email: {\tt tobias.jager@college-de-france.fr}}}

\newcommand{\homeo}{\ensuremath{\mathrm{Homeo}}}
\newcommand{\TT}{\ensuremath{\mathbb{T}}}

\begin{document}

\setlength{\oddsidemargin}{-0.03\textwidth}
\setlength{\evensidemargin}{-0.03\textwidth}

\maketitle 

\abstract{A conservative irrational pseudo-rotation of the two-torus is
  semi-conjugate to the irrational rotation if and only if it has the property
  of {\em bounded mean motion} \cite{jaeger:2008a}. (Here `irrational
  pseudo-rotation' means a toral homeomorphism with a unique and totally
  irrational rotation vector.) The aim of this note is to explore this concept
  some further. For instance, we provide an example which shows that the
  preceding statement does not hold in the non-conservative case. Further, we
  collect a number of observations concerning the case where the bounded mean
  motion property fails. In particular, we show that a non-wandering irrational
  pseudo-rotation of the two-torus with unbounded mean motion has sensitive
  dependence on initial conditions.}

\section{Introduction}
\label{Intro}

We denote by $\homeo_0(\TT^d)$ the set of homeomorphisms of the $d$-dimensional
torus which are homotopic to the identity. An important topological and
dynamical invariant for this class of maps is given by the rotation set. Given
$f \in \homeo_0(\TT^d)$ with lift $F : \R^d \to \R^d$, it is defined as
\begin{equation}\label{eq:1}
  \rho(F) \ := \ \left\{ \rho \in \R^d \left|\  \exists n_i\nearrow\infty,\ z_i\in\R^d :
 \iLim (F^{n_i}(z_i) - z_i)/n_i \ = \ \rho \right.\right\} \ .
\end{equation}
The set $\rho(F)$ is always compact, and if $d=2$ it also convex
\cite{misiurewicz/ziemian:1989}. At least in dimension two, it is well-known
that the shape of $\rho(F)$ is closely related to the dynamical properties of
$f$. For instance, if $\rho(F)$ has non-empty interior, then all rational
rotation vectors in $\inte(\rho(F))$ are realised by periodic orbits
\cite{franks:1989}, the topological entropy of $f$ is strictly positive
\cite{llibre/mackay:1991} and the rotation set is continuous in $f$ with respect
to the topology of uniform convergence \cite{misiurewicz/ziemian:1991} (in
general, the dependence on $f$ is only upper semi-continuous). 

In contrast to this, the situation is less well-understood when the rotation set
is {\em `thin'}, meaning that it has empty interior. A particular case is that
of an {\em (irrational) pseudo-rotation}, which is defined as a toral
homeomorphism $f$ whose rotation set is reduced to a single (totally irrational)
rotation vector $\rho$. In this situation, we say $f$ has {\em bounded mean
  motion}, if the {\em deviations from the constant rotation}
\begin{equation}
  \label{eq:2}
D(n,z) \ := \ F^n(z) - z - n\rho \ 
\end{equation}
are uniformly bounded in $\theta,x$ and $n$. This notion already played an
important role in the description of the dynamics of certain skew product
transformations, like quasiperiodically forced circle homeomorphisms
\cite{stark/feudel/glendinning/pikovsky:2002,jaeger/stark:2006,bjerkloev/jaeger:2006a}
and almost periodically forced circle flows \cite{yi:2004,huang/yi:2007}. The
fact that it is equally important for understanding, and classifying, the
possible dynamics of irrational pseudo-rotations is documented by the following
statement.
\begin{thm}[\cite{jaeger:2008a}] \label{t.linearization} Suppose $f\in\homd$ is
  minimal, or $d=2$ and $f$ is conservative. Then $f$ is semi-conjugate to an
  irrational rotation on the $d$-dimensional torus if and only if it is an
  irrational pseudo-rotation with bounded mean motion.
\end{thm}
At least of the minimal/conservative case, this result suggests to divide
irrational pseudo-rotations into two classes, those with and those without
bounded mean motion. However, the usefulness of such a dichotomy will strongly
depend on the information which can be deduced from the unboundedness
of mean motion. Furthermore, it would certainly be desireable to weaken the
recurrence assumption if possible.

In other words, Theorem~\ref{t.linearization} raises two obvious questions: (i)
What happens in the non-minimal/dissipative case? (ii) What happens if the
bounded mean motion property fails? While being far from giving a complete
picture, the aim of this note is to collect a number of results and observations
related to these two questions.

Concerning the first, we will provide examples of irrational
pseudo-rotations with bounded mean motion which are not semi-conjugate to an
irrational rotation (see Section~\ref{BoundedMeanMotion}). These examples have a
very simple structure, namely they are skew-products over Denjoy
counterexamples, with rotations on the fibres. Hence, they are still
semi-conjugate to a one-dimensional irrational rotation, and they have very bad
recurrence properties, since they admit homotopically non-trivial wandering open
sets. We have to leave open here whether more sophisticated examples without
these two properties exist. In particular, it is still possible that the
assumptions in Theorem~\ref{t.linearization} could be weakened to some extent,
for example to the non-existence of wandering open sets.

In order to address the second question, we will first of all consider the
behaviour of the quantities $D_n(\theta,x)$ themselves in the case of unbounded
mean motion. This is motivated by the fact that a good understanding of this
behaviour already turned out to be crucial in the theory of the skew product
transformations mentioned above. However, in order to obtain a good picture,
some recurrence assumption on the system is needed. The best description can be
given for minimal systems (see Section~\ref{UnboundedMeanMotion}). Skew-products
over Denjoy counterexamples allow again to demonstrate that a number of natural
properties that hold in the minimal case cannot be expected in general.

Finally, in Section \ref{LyapunovStability} we show that for an irrational
pseudo-rotation of the two-torus without wandering open sets, the unboundedness
of mean motion implies sensitive dependence on initial conditions. The interest
of this result lies in particular in the fact that it demonstrates how the
unboundedness of mean motion allows to obtain further information about the
dynamics of the system. \medskip

\noindent {\bf Notation and terminology.}
For the sake of readability, we will use a number of conventions which might not
be absolutely consistent from a strictly formal point of view, but should not
cause any ambiguities and greatly simplify notation. For example, we will
identify points in $\T^d$ with their lifts whenever the particular choice of the
lift does not matter. Whenever there is a canonical distance function on a
metric space, we will denote it by $d$. In the case of $\R,\ \R^d$ or $\T^d$,
this will always be the usual Euclidean distance. The $\eps$-neighbourhood of a
point $z$ in a metric space $X$ will be denoted by $B_\eps(z) = \{x\in X\mid
d(z,x) < \eps\}$. The Euclidean length of a vector $v\in\R^d$ will simply be
denoted by $\|v\|$. Given a vector $v\in\R^d$ and a subset $S\ssq \R^d$, we let
$S+v := \{s+v\mid s\in S\}$.

Quotient maps like $\R \to \T^1$ or $\R^d \to \T^d$ will all be denoted by
$\pi$, in product spaces $\pi_i$ will be the projection to the respective
coordinate.  The rotation on $\T^d$ by an element $\rho\in\T^d$ will be denoted
by $R_\rho$, as well as the translation by a vector $\rho\in \R^d$ on $\R^d$.

\medskip

\noindent {\bf Acknowledgements.} I would like to thank Sylvain Crovisier and
Fr\'ed\'eric Le Roux for stimulating discussions. Further, I would like to thank
J.-C.\ Yoccoz and the Coll\`ege de France for their hospitality and support
during a long-term visit, and acknowledge support from the German Science
Foundation (research fellowship Ja 1721/1-1).


\section{Counter-examples to the existence of a semi-conjugacy}
\label{BoundedMeanMotion}

In this section, we provide examples of irrational pseudo-rotations which have
bounded mean motion, but are not semi-conjugate to an irrational rotation. The
main idea is to use skew products over Denjoy counter-examples, with rotations
on the fibres. The dynamics over the wandering intervals can then be modified,
without affecting the rotation set, in order to produce some irregular
behaviour. Since we will employ this basic construction several times, we first
want to give a brief outline.

Given a rotation vector $\rho = (\rho_1,\rho_2) \in (\R\smin\Q) \times \R$, we choose
an orientation-preserving circle homeomorphism $g$ with rotation number
$\rho(g)=\rho_1$ and unique minimal set $M \subsetneq \kreis$. Further, we choose
some continuous function $\beta:\kreis \to \R$ which satisfies $\beta(x) = 0
\ \forall x \in M$ and define
\begin{equation}
  \label{eq:3}
f : \kreis \times \R \to \kreis \times \R \quad , \quad (x,y) \mapsto
(g(x),y+\rho_2+\beta(x) \bmod 1) \ .
\end{equation}
This yields an irrational pseudo-rotation $f$ with unique rotation vector
$\rho$. One way to see this is to note that any $f$-invariant measure $\mu$ must
project down to the unique $g$-invariant measure $\nu$. Hence, the topological
support of $\mu$ is contained in the set $M \times \kreis$. On this set, the
replacement function $D(1,z) = F(z)-z$ (defined via a suitable lift $F$ of $f$)
is constant and equal to $\rho$, and therefore $\int_{\ntorus} F(z)-z \ d\mu =
\rho$ for all $f$-invariant measures $\mu$. This implies $\rho(F) = \{\rho\}$ as
claimed. The important point is that on the set $M^c$ the function $\beta$ can
still be chosen arbitrarily in order to produce the desired phenomena.

As a first application of this construction, we prove the following:
\begin{prop}
  \label{p.semiconj-counterex} Given any totally irrational rotation vector
  $\rho=(\rho_1,\ldots,\rho_d) \in \R^d, \ d\geq 2$, there exists an irrational
  pseudo-rotation $f \in \homd$ which has rotation vector $\rho$ and bounded
  mean motion, but which is not semi-conjugate to the irrational rotation
  $R_\rho$.
\end{prop}
We will actually carry out the construction only in dimension $d=2$, but the
modifications for the general case are minor. Before we turn to the proof, we
need the following lemma.
\begin{lem}
  \label{l.fibred-conjugacy}
  Let $\rho = (\rho_1,\rho_2) \in \R^2$ be totally irrational. Suppose $g \in
  \homeo_0(\T^1)$ is semi-conjugate to $R_{\rho_1} : \kreis \to \kreis$ by a
  semi-conjugacy $\phi$ and $f\in \homeo_0(\T^2)$ is semi-conjugate to
  $R_\rho:\torus\to\torus$ by a semi-conjugacy $\psi$. Further, assume that $f $
  is a skew-product over $g$, that is $\pi_1\circ f(x,y) = g(x) \ \forall (x,y)
  \in \T^2$.
\alphlist
\item If $f$ has a unique minimal set, then the semi-conjugacy $\psi$ can be chosen
  such that $\pi_1\circ \psi(x,y) = \phi(x) \ \forall (x,y) \in \T^2$. 
\item If $f$ is given by (\ref{eq:3}), with $\beta_{|M} = 0$, then $\psi$ can be
  chosen such that $\psi(x,y) = (\phi(x),y) \ \forall (x,y) \in M \times \T^1$.
  \listend
\end{lem}
\proof The argument is a slight variation of the one in \cite[Section
4.14]{herman:1983}. Choose lifts $G,\Phi,F$ and $\Psi$ of the maps $g,\phi,f$
and $\psi$ to $\R$ and $\R^2$, respectively, such that $\Phi \circ G =
R_{\rho_1} \circ \Phi$, $\Psi\circ F = R_\rho \circ \Psi$ and $\pi_1\circ F(x,y)
= G(x)$. We claim that the function $\eta = \pi_1 \circ \Psi - \Phi \circ \pi_1$
is constant. In order to see this, note first that
\begin{eqnarray*}
  \eta \circ F(x,y) & = & \pi_1 \circ \Psi \circ F(x,y) - \Phi \circ \pi_1 \circ F(x,y) \\
  &  = & \pi_1 \circ R_\rho \circ \Psi(x,y) - \Phi \circ G(x,y) \\
  & = & R_{\rho_1} \circ \pi_1 \circ \Psi(x,y) - R_{\rho_1} \circ \Phi(x) \\
  & = & \pi_1 \circ \Psi(x,y) - \Phi(x) \ = \ \eta(x,y) \ .
\end{eqnarray*}
Hence $\eta$ is $F$ invariant, and we can interpret it as a $f$-invariant
function $\T^2 \to \R$. Since $\eta$ is continuous, it follows first that it
must be constant on the unique minimal set of $f$, and since any point must
contain this minimal set in its $\omega$-limit, $\eta$ must be constant
everywhere. Composing $\Psi$ with a rotation in the $x$-direction, we can assume
that $\eta=0$. This proves part (a).

Now assume that $f$ is of the form (\ref{eq:3}). Let $\mu = \Leb_{\T^1} \circ
\phi$. Then $\mu$ is $g$-invariant, and $\mu \times \Leb_{\kreis}$ is
$f$-invariant. Furthermore, $\psi$ has to map $\mu \times \Leb_{\kreis}$ to the
unique $R_\rho$-invariant measure $\Leb_{\T^2}$. This means that for $\mu$-a.e.\
$x$, the map $\psi_x : y \mapsto \pi_2\circ \psi(x,y)$ must preserve Lebesgue
measure and hence be a rotation. By continuity, this is true for all $x\in
\supp(\mu) = M$. Using the notation $f_x(y) = \pi_2\circ f(x,y)$, this further
implies that for all $x \in M$ there holds
\begin{eqnarray*}
  \psi_{g(x)}(y) + \rho_2  & = & \psi_{g(x)}(y+\rho_2) \ = \ \psi_{g(x)} \circ f_x(y) \\
  & = &  R_{\rho_2} \circ \psi_x(y) \ = \ \psi_x(y) + \rho_2 \ .
\end{eqnarray*}
Thus the rotations $\psi_x$ and $\psi_{g(x)}$ are the same, and using continuity
of $\psi$ and minimality of $M$ once more we obtain that $x \mapsto \psi_x$ is
constant on $M$. Composing $\psi$ with a rotation in the $y$-direction, we may
therefore assume that $\psi_x = \Id_{\kreis} \ \forall x \in M$, which proves
statement~(b).

\qed

\proof[Proof of Proposition~\ref{p.semiconj-counterex}] As mentioned, we give
the proof in dimension 2.  Let $\rho = (\rho_1,\rho_2)$ be totally irrational
and define $f$ as in (\ref{eq:3}), where we choose $\beta$ as follows. First,
let $\nfolge{\alpha_n}$ be a sequence of real numbers, such that
\romanlist
\item $\nLim \alpha_n = 0$;
\item $\sup_{n\in\Z} \left| \sum_{i=0}^{n} \alpha_i \right| \leq 1$;
\item $\limsup_{n\to\infty} \sum_{i=0}^n \alpha_i = 1$;
\item $\liminf_{n\to\infty} \sum_{i=0}^n \alpha_i = 0$.  \listend
Let $I_0=[a_0,b_0]$ be a connected component of $M^c$ and let $I_n = [a_n,b_n] =
g^n(I_0)$. Further, fix some $x_0 \in (a_0,b_0)$ and let $x_n =
g^n(x_0)$. Choose a continuous function $\beta_0 : I_0 \to [0,1]$ which
satisfies $\beta_0(a_0)=\beta_0(b_0)=0$ and $\beta_0(x_0) = 1$. Then define
\[
\beta(x) \ := \ \left\{ \begin{array}{cc}
    \alpha_n \beta_0(g^{-n}(x)) & \mathrm{if}\   x \in I_n \textrm{ for some } n\in\Z \\
    0 & \mathrm{otherwise}\end{array} \right. \quad .
\]
Property (ii) of the $\alpha_n$ then implies that $f$ has bounded mean
motion. We claim that $f$ cannot be semi-conjugate to the irrational rotation
$R_\rho$. Suppose for a contradiction that $\psi$ is a semi-conjugacy between
$f$ and $R_\rho$. By Lemma~\ref{l.fibred-conjugacy}, we can assume that $\psi$
is of the form $\psi(x,y) = (\phi(x),\psi_x(y))$, where $\phi$ is a
semi-conjugacy between $g$ and $R_{\rho_1}$, and $\psi_x = \Id_{\kreis} \
\forall x \in M$. Due to uniform continuity, there exists $\delta > 0$, such
that for all $z,z' \in \torus$ with $d(z,z') < \delta$ there holds
$d(\psi(z),\psi(z')) < \achtel$. If we choose $k\in \N$ such that $d(x_k,b_k) <
\delta$ (note that $|I_k| \to 0$ as $k \to \infty$), then it follows that for
all $n\in\Z$ there holds
\begin{eqnarray}
  \lefteqn{d(\psi \circ f^n(x_k,0),(\phi(b_{k+n}),n\rho_2)) \ = }  \nonumber \\
  & = & d(\psi \circ f^n(x_k,0),\psi \circ f^n(b_k,0)) \label{eq:4}
  \\ & = & d\left(R_\rho^n\circ \psi(x_k,0),R_\rho^n\circ\psi(b_k,0)\right) \
  < \ 1/8 \ . \nonumber  
\end{eqnarray}
On the other hand, there holds $\pi_2 \circ f^n(x_k,0) = \sum_{i=k}^{k+n-1}
\alpha_i \bmod 1$. Hence, due to properties (i), (iii) and (iv) above there
exist infinitely many $n\in\N$, such that
\[
d\left(\pi_2\circ f^n(x_k,0),n\rho_2+1/2\right) \ < \ \delta/2 \ .
\]
If we choose such an $n$ with the additional property that $d(x_{k+n},b_{k+n}) <
\frac{\delta}{2}$, then
\[
d(f^n(x_k,0),(b_{k+n},n\rho_2+1/2)) \ < \ \delta \ .
\]
Since $\psi(b_{k+n},n\rho_2+1/2) = (\phi(b_{k+n}),n\rho_2+1/2)$, it follows that
\begin{eqnarray}
d(\psi\circ f^n(x_k,0),(\phi(b_{k+n}),n\rho_2+1/2)) \ < \ 1/8 \ .
\end{eqnarray}
This clearly contradicts (\ref{eq:4}). 

\qed
\medskip


\section{Unbounded Mean Motion}
\label{UnboundedMeanMotion}

We first want to mention that the existence of irrational pseudo-rotations with
unbounded mean motion is well-known. In \cite{herman:1983}, Herman provides
examples of homeomorphisms of the two-torus which are skew products over an
irrational rotation (so-called quasiperiodically forced circle homeomorphisms)
and have totally irrational rotation vector, but are not semi-conjugate to the
corresponding irrational rotation. Since for such skew products the statement of
Theorem~\ref{t.linearization} holds without any recurrence assumption
\cite{jaeger/stark:2006}, these examples have unbounded mean motion. Other
examples of irrational pseudo-rotations with unbounded mean motion are provided by uniquely
ergodic and weakly mixing torus diffeomorphisms, as constructed for example in
\cite{fayad:2002}.

In order to give a more detailed discussion, we first need a refined version of
the notions introduced in Section~\ref{Intro}. Given $f\in\homd$ with lift $F$
and a subset $A \ssq \T^d$, we define the rotation set of $f$ on $A$ as
\begin{equation}
  \label{eq:5}
  \rho_A(F) \ := \ \left\{ \rho \in \R^d \left|\ \exists n_i \nearrow \infty,\ x_i
  \in A : \iLim (F^{n_i}(x_i) - x_i)/n_i = \rho \right. \right\} \ .
\end{equation}
Further, when $v$ is a vector in $\R^d\smin\{0\}$, $\lambda \in \R$ and
$\rho_A(F) \ssq \lambda v+\{v\}^\perp$, then for any $z\in A$ we define the {\em
  deviations from the constant rotation parallel to $v$} as
\begin{equation}
  \label{eq:2a}
  D_v(n,z) \ := \ \langle F^n(z)-z-n\rho,v\rangle \ .
\end{equation}
Here $\rho$ can be any vector in $\lambda v + \{v\}^\perp$. When the quantities
$D_v(n,z)$ are uniformly bounded, we say $f$ has {\em bounded deviations
  parallel to $v$ on $A$}. 

A priori, it is not entirely clear whether there exist irrational
pseudo-rotations with unbounded mean motion in all directions (that is, parallel
to all vectors $v\in\R^d\smin\{0\}$. For the examples mentioned above, there
always exists some vector $v$, such that the deviations parallel to $v$ are
bounded.  Hence, for the sake of completeness, we will provide examples with
unbounded deviations parallel to all vectors at the end of this section (see
Proposition~\ref{p.unbounded-examples}).

Before, we want to gather some more information about the behaviour of
the quantities $D_v(n,z)$ when the mean motion parallel to $v$ is not
bounded. In spirit, the results we present below are close to those collected in
\cite{stark/feudel/glendinning/pikovsky:2002} for quasiperiodically forced
circle homeomorphisms. 
\medskip

In the whole section, we suppose that $f\in\homd$, $F$ is a lift of $f$ and
$K\ssq \T^d$ is a compact $f$-invariant subset. Further, we always assume that
$\lambda$ is a real number and $v$ is a vector in $\R^d\smin\{0\}$. As a first
observation, we note that there always exists orbits on which the deviations are
semi-bounded in the following sense.
\begin{lem}
\label{l.semi-bounded-orbits}
Suppose that $\rho_K(F) \ssq \lambda v + \{v\}$. Then given any $\eps>0$, there
exist $z_-,z_+ \in K$, such that
\[
\sup_{n\in\N} D_v(n,z_-) < \eps \quad \textrm{and} \quad \inf_{n\in\N}
D_v(n,z_+)>-\eps \ .
\]
\end{lem}
\proof Without loss of generality, we may assume $\|v\|=1$. Suppose for a
contradiction that for all $z\in K$ there holds
\[
\sup_{n\in\N} D_v(n,z) \ \geq \ \eps \ > \ 0 \ .
\]
Then for any $z\in K$, there exists an integer $n(z)$ and a small
neighbourhood $U(z)$, such that
\[
\langle F^{n(z)}(z')-z',v\rangle -n(z)\lambda \ \geq \ \eps/2 \quad \forall z'
\in U(z) \ .
\]
Using compactness, we can find $z_1 \ld z_m$ with $K \ssq \bigcup_{i=1}^m
U(z_i)$. Let $M := \max_{i=1}^m n(z_i)$. Inductively define a sequence of integers
$N_j$ and points $\zeta_j$ as follows: Let $N_0=0$ and $\zeta_0$ be arbitrary. If
$N_0\ld N_j$ and $\zeta_0\ld \zeta_j$ are given, choose $k\in\{1\ld m\}$ with $\zeta_j \in
U(z_k)$ and let $\zeta_{j+1} := F^{n(z_k)}(\zeta_j)$ and $N_{j+1} := N_j + n(z_k)$. Then
\begin{eqnarray*}
  \lefteqn{\left\langle (F^{N_j}(z_0)-z_0),v\right\rangle  \ = \ \langle z_j-z_0,v\rangle} \\ & = & 
  \sum_{i=0}^{j-1} \langle z_{i+1}-z_i,v\rangle \ \geq \  N_j\lambda + (j\eps)/2 \ \geq \
  (\lambda+\eps/2M) \cdot  N_j \  .
\end{eqnarray*}
By going over to a subsequence if necessary, we may further assume that
$(F^{N_j}(z_0) - z_0)/N_j \to \rho \in \rho(F)$ as $j\to \infty$. Thus we obtain
\begin{eqnarray*}
  \langle \rho,v\rangle & = & \jLim \langle (F^{N_j}(z_0)-z_0)/N_j,v\rangle \ \geq \ \lambda+\eps/2M \ ,
\end{eqnarray*}
contradicting the assumptions. 

\qed
\medskip

In general, it seems difficult to say much more about the qualitative behaviour
of the quantities $D_v(n,z)$. For instance, even when the deviations from the
constant rotation are not uniformly bounded, they may be bounded on every single
orbit, as the following example shows.
\begin{example} \label{e.unbounded} \alphlist
\item
  Suppose $f\in\homeo_0(\torus)$ is a skew product over a Denjoy counter\-example
  $g$, of the same form as in (\ref{eq:3}). Let $I=(a,b)$ be a connected component
  of $M^c$, where $M$ is the unique minimal set of $g$, and $I_n = g^n(I)$. In
  order to define $\beta$ in (\ref{eq:3}), choose a continuous function $\beta_0
  : [a,b] \to [0,1]$ with $\beta_0(a)=\beta_0(b)=0$ and $\beta(x_0) = 1$ for some
  $x_0\in I$. Let
  \[
      \beta(x) \ := \ \left\{
        \begin{array}{cc}
          \beta_0\circ g^{-n}(x)/|n| & \textrm{if }\   x \in I_n \textrm{ with } n\leq 0 \\
          0 & otherwise
        \end{array} \ . \right.
  \]
  Then it is easy to see that $f$ does not admit bounded mean motion, but
  $\sup_{n\in\N}\|D(n,z)\| < \infty \ \forall z\in\torus$.
\item
  Evidently, the preceeding example still exhibits orbits with unbounded mean
  motion backwards in time, that is $sup_{n\in\N}\|D(-n,z)\| = \infty$ for some
  $z$. However, even this need not be true, as a slight modification of the
  construction shows. Suppose $I^k,\ k\in\N$ is a sequence of disjoint open
  subintervals of $I$ and let $I^k_n := g^n(I^k)$. Choose $\beta_0: I \to
  [0,1]$, such that $\sup_{x\in I^k} \beta_0(x) = 1/k$. Then define 
  \[
      \beta(x) \ := \ \left\{
        \begin{array}{cc}
          \beta_0\circ g^n(x)/n & \textrm{if }\   x \in I^k_n \textrm{ with } |n| \leq 2^{k^2} \\
          0 & otherwise
        \end{array} \ . \right.
  \]
  Since $\sum_{n=1}^{2^{k^2}} 1/n \geq k^2$, we have $\sup_{n\in\N} \|D(n,(x,0))\|
  \geq k$ for some $x\in I^k$, and therefore $\sup_{n\in\N,z\in\torus}\|D(n,z)\| =
  \infty$, but on the other hand there holds $\sup_{n\in\Z}\|D(n,z)\| < \infty \
  \forall z\in\torus$.
\item A slight modification of these constructions also shows that the
  deviations from the constant direction can be unbounded in one direction, but
  at the same time bounded in the opposite direction. More precisely, it is
  possible to have $\sup_{n\in\N} D_v(n,z) = \infty$ for some $z\in\T^d$, but
  $\inf_{n\in\Z,z\in\T^d} D_v(n,z) > -\infty$, and vice versa.
\listend
\end{example}
The situation becomes different when $f$ has some topological recurrence
properties. Starting point is the following simple observation. Let ${\cal
  O}^+(z) = \{ f^n(z) \mid n\in \N\}$ denote the forward orbit of $z$.
\begin{lem}
  \label{l.orbitclosure} Suppose $\rho(F) \ssq \lambda v + \{v\}^\perp$ and
  $\sup_{n\in\N} |D_v(n,z)| \leq C < \infty$ for some $z \in \T^d$. Then
  $\sup_{n\in\N} |D_v(n,z)| \leq 2 C \ \forall z \in \overline{{\cal O}^+(z)}$.
\end{lem}
\proof Suppose $D_v(n,z_0) > 2C$ for some $z_0 \in \overline{{\cal
    O}^+(z)}$. Then for some sufficiently small neighbourhood $U$ of $z_0$,
there holds $D_v(n,z') > 2C \ \forall z'\in U$. In particular, when $f^k(z) \in
U$, then $D_v(n,f^k(z)) > 2C$. However, this contradicts
\[
  D_v(n,f^k(z)) \ = \ D_v(n+k,z) - D_v(k,z) \ \leq \ 2C \ .
\]
\qed
\medskip

Now suppose $f_{|K}$ is topologically transitive and denote the set of all
points in $K$ with dense forward orbit by ${\cal T}_K$. Recall that ${\cal T}_K$
is a residual (meaning dense $G_\delta$) subset of $K$.
\begin{cor}
  \label{c.transitive-orbits} Suppose $f_{|K}$ is topologically transitive and
  $\rho_K(F) \ssq \lambda v + \{v\}^\perp$.  If there holds $\sup_{n\in\N,z\in K}|D_v(n,z)|
  = \infty$, then $\sup_{n\in\N}|D_v(n,z)| = \infty \ \forall z\in{\cal T}_K$.
\end{cor}
This statement is particularly interesting, when $f_{|K}$ is minimal, since this
means that ${\cal T}_K=K$ and therefore all orbits have unbounded
deviations. Together with Lemma~\ref{l.semi-bounded-orbits}, this entails
the following corollary.
\begin{cor}
  \label{c.sensitive-dependence}
  Suppose $f$ is minimal, $\rho_K(F) \ssq \lambda v+\{v\}^\perp$ and the deviations from
  the constant rotation parallel to $v$ are not uniformly bounded on $K$. Then
  $f_{|K}$ has sensitive dependence on initial conditions.
\end{cor}
\proof Due to uniform continuity, there exists $\eps > 0$, such that $d(z_1,z_2)
< \eps$ implies $d(F(z_1),F(z_2)) < 1/2$. Suppose $\delta > 0$ and $z\in\T^d$
are given. We have to show that there exist $z_1,z_2 \in B_\delta(z)$ and $n\in
\N$, such that $d(f^n(z_1),f^n(z_2)) \geq \eps$. In order to do so, choose
$z_1,z_2\in B_\delta(z)$ with the properties that $\sup_{n\in\N} D_v(n,z_1) <
\infty$ and $\inf_{n\in\N} D_v(n,z_2) > -\infty$. This is possible due to
Lemma~\ref{l.semi-bounded-orbits}, minimality, and the fact that these
properties are invariant under iteration by $f$.

Since the deviations are unbounded on all orbits by
Lemma~\ref{c.transitive-orbits}, the distance between $F^n(z_1)$ and $F^n(z_2)$
in $\R^d$ must become arbitrary large. The choice of $\eps$ above therefore
implies that for some $n\in\N$ there holds $d(F^n(z_1),F^n(z_2)) \in
[\eps,1/2)$. However, as long as two points have distance less than $1/2$ in
$\R^d$, their projections have exactly the same distance in $\T^d$. Hence
$d(f^n(z_1),f^n(z_2)) \geq \eps$ as required.

\qed
\medskip

The preceding statement provides a first indication of how information on the
behaviour of the deviations $D(n,z)$ may be used in order to draw further
conclusions about the dynamical properties of $f$. In the proof, we have used
that on the one hand the vectors $D(n,z)$ become arbitrarily large, but on the
other hand we have also exploited the fact that they grow in different
directions, depending on the starting point, due to
Lemma~\ref{l.semi-bounded-orbits}.

In other words, we have not only made use of information on the size, but also
on the behaviour of the directions of the vectors $D(n,z)$. In the following, we
want to obtain a more detailed description of this behaviour. It turns out that
this requires a much more careful analysis, and we will here restrict to the case
where the dynamics are minimal. Our aim is the following statement.

\begin{thm} 
\label{t.unbounded-deviations}
Suppose that $f_{|K}$ is minimal and there exists a linear subspace $V\ssq \R^d$,
such that $\rho_K(F) \ssq v_0 + V^\perp$ for some vector $v_0\in V$. Further,
assume that for any $v\in V$ the deviations parallel to $v$ are unbounded on
$K$. Then there exists a residual set ${\cal R} \ssq K$, such that for all
$z\in{\cal R}$ there holds
\[
\sup_{n\in\N} D_v(n,z) \ = \ -\inf_{n\in\N} D_v(n,z) \quad \forall v \in V \ .
\]
\end{thm}
Needless to say that for a minimal irrational pseudo-rotation with unbounded
deviations in all directions, as we will construct it below, this implies a
rather `wild' behaviour. Our first step is to consider a one-dimensional affine
subspace.
\begin{lem}
  \label{l.minimal-fully-unbounded} Suppose $f_{|K}$ is minimal, $\rho_K(F) \ssq
  \lambda v + \{v\}^\perp$ and the deviations from the constant rotation parallel to $v$
  are not uniformly bounded. Then there exists a residual subset ${\cal R} \ssq
  K$, such that
  \[
  \sup_{n\in\N} D_v(n,z) \ = \ -\inf_{n\in\N}D_v(n,z) \ = \ \infty \quad \forall
  z\in{\cal R} \ .
  \]
\end{lem}
\proof Note that there holds
\[ {\cal R}^+ := \left\{ z \in\T^d \left|\ \sup_{n\in\N} D_v(n,z) = \infty \right.\right\} \ = \
\bigcap_{k\in\N} {\cal R}^+_k \ ,
\]
where ${\cal R}^+_k := \left\{ z \in \T^d \mid \exists n\in\N : D_v(n,z) >
  k\right\}$. The sets ${\cal R}^+_k$ are obviously open, and from
Lemma~\ref{l.semi-bounded-orbits} and the minimality of $f$ it follows that they
are also dense (they contain the orbit of the point $z_+$ in
Lemma~\ref{l.semi-bounded-orbits}). Hence ${\cal R}^+$ is residual. The same
applies to the set ${\cal R}^- := \left\{z\in\T^d \left|\ \inf D_v(n,z) = -\infty
  \right. \right\}$, and the intersection yields the required set ${\cal R}$.

\qed
\medskip

In order to generalise Lemma~\ref{l.minimal-fully-unbounded} to
higher-dimensional hyperplanes, it will be convenient not only to consider the
scalar quantities $D_v(n,z)$, but also the directions of the vectors $D(n,z)$
and their projections to some linear subspace $V \ssq R^d$. If $\pi_V : \R^d \to
V$ denotes the orthogonal projection to $V$, let
\begin{equation}
D_V(n,z) \ := \ \pi_V(D(n,z)) \ .
\end{equation}
As usual, the unit sphere in $\R^d$ is denoted by $\sph^{d-1} := \{ v\in\R^d
\mid \|v\| = 1\}$. Given any $v\in\sph^{d-1}$ we let $S^+(v) := \{
w\in\sph^{d-1} \mid \langle v,w\rangle > 0\}$. We define the set of {\em
  `asymptotic directions'} of $z\in K$ as
\[
  \label{eq:10} {\cal A}_V(z) \ := \ \left\{ v\in\sph^{d-1} \left|\
  \forall \eps,r> 0 \
  \exists n\in\N: \|D_V(n,z)\| > r \textrm{ and } \frac{D_V(n,z)}{\|D_V(n,z)\|} 
 \in B_\eps(v)\right. \right\} \ ,
\]
and let ${\cal A}^K_V := \bigcup_{z\in K} {\cal A}_V(z)$ and ${\cal R}_V(v) :=
\{ z\in\T^d \mid v \in {\cal A}_V(z)\}$. The following lemma collects some
elementary properties of these sets.
\begin{lem}
  \label{l.asymptotic-directions}~ Suppose $V \ssq \R^d$ is a linear subspace
  and $\rho_K(F) \ssq v_0 + V^\perp$ for some $v_0\in V$. Then the following holds:
\alphlist
\item For all $z\in K$, the set ${\cal A}_V(z)$ is compact and ${\cal A}_V(f(z)) =
  {\cal A}_V(z)$. 
\item If $\sup_{n\in\N} \|D_V(n,z)\| = \infty$, then
  ${\cal A}_V(z)$ is non-empty.
\item If $\sup_{n\in\N}D_v(n,z) = \infty$ for some $v\in V$, then
  $\overline{S^+(v)} \cap {\cal A}_V(z) \neq \emptyset$.
\item ${\cal R}_V(v)$ is an invariant $G_\delta$-subset of $K$. In particular, if
  ${\cal R}_V(v)$ contains a point with dense forward orbit, then it is a residual
  set.  \listend
\end{lem}
\proof \alphlist
\item The fact that ${\cal A}_V(f(z)) = {\cal A}_V(z)$ is obvious from the
  definition. In order to see that ${\cal A}_V(z)$ is compact, let
\[
   A_k(z) \ := \ \left\{ v\in\sph^{d-1} \left| \exists n\in\N : \|D(n,z)\| > k
    \textrm{ and } \frac{D(n,z)}{\|D(n,z)\|} \in B_{1/k}(v) \ \right.\right\} \ . 
\]
Then ${\cal A}_V(z) = \kcap \overline{A_k(z)}$ is compact as an
intersection of compact sets. 
\item When the sequence $\nfolge{D_V(n,z)}$ is unbounded,
then $\overline{A_k(z)} \neq \emptyset \ \forall k\in\N$, and hence the intersection
is non-empty.
\item The fact that $\sup_{n\in\N} D_v(n,z) = \infty$ implies that for all
  $k\in\N$, the sets $\overline{S^+(v)} \cap\overline{A_k(z)}$ are non-empty, and again
  this carries over to their intersection.
\item The invariance of ${\cal R}_V(v)$ follows directly from ${\cal A}_V(f(z)) =
  {\cal A}_V(z)$. In order to see that ${\cal R}_V(v)$ is a $G_\delta$-set, let
  $R_k(v) := \{ z\in\T^d \mid v \in A_k(z)\}$. It is easy to see that these sets
  are open, and ${\cal R}_V(v) = \kcap R_k(v)$. 
\listend

\qed
\medskip

Further, we need a piece of linear algebra.
\begin{lem} \label{l.hyperplanes}
  Suppose $V\ssq \R^d$ is a linear subspace of dimension $k\leq d$ and $A \ssq
  \sph^{d-1} \cap V$ is a compact set with the property that 
  \[
  \sph^{d-1} \cap V \ \ssq \ \bigcup_{v\in A} \overline{S^+(v)} \ .
  \]
  Then there exist vectors $v^1 \ld v^n \in A$ and a linear subspace $\tilde V
  \ssq V$ of dimension $< k$, such that there holds
  \[
  \left(\sph^{d-1} \cap V\right) \smin \incup S^+(v_i) \ \ssq \ \tilde V \ .
  \]
\end{lem}
\proof We proceed by induction on the dimension $k$ of $V$. If $k=1$, such that
$\sph^{d-1} \cap V$ consists of exactly two antipodal points, then the statement
is obvious. Hence, suppose that it holds for all linear subspaces of of
dimension $< k$. If
\[
  \sph^{d-1} \cap V \ \ssq \ \bigcup_{v\in A} S^+(v) \ ,
\]
then using compactness we can find $v^1 \ld v^n \in A$, such that 
\[
  \sph^{d-1} \cap V \ \ssq \ \incup S^+(v^i)  
\]
and we are finished. (In this case $\dim \tilde V = 0$.) 

Otherwise, there exists $w\in \sph^{d-1} \cap V$, such that $S^+(w) \cap A =
\emptyset$. In this case, let $V_0 := V \cap \{w\}^\perp$ and $A_0 := A \cap
V_0$. Then there holds
\[
    \sph^{d-1} \cap V_0 \ \ssq \ \bigcup_{v\in A_0} \overline{S^+(v)} \ ,
\]
since for any $u \in  \sph^{d-1} \cap V_0$ the set 
\[
\overline{S^+(u)} \cap A_0 \ = \ \bigcap_{\eps > 0}
\overline{S^+\left(w_\eps\right)} \cap A \ ,
\]
where $w_\eps = \frac{w+\eps u}{\|w+\eps u\|}$, is non-empty as the intersection
of a nested sequence of non-empty compact sets. Applying the induction
assumption to $V_0$ yields vectors $v^1\ld v^n \in A_0$ and a linear subspace
$\tilde V_0 \ssq V_0$ of dimension $<k-1$, such that
\[
\left(\sph^{d-1} \cap V_0\right)\setminus \bigcup_{i=1}^n S^+(v^i) \ \ssq \  \tilde V_0 \ .
\]
If we now define $\tilde V$ as the linear subspace spanned by $\tilde V_0$ and
$w$, then the assertions of the lemma are satisfied.

\qed
\medskip

\proof[\bf Proof of Theorem~\ref{t.unbounded-deviations}] We proceed by induction on
the dimension $k$ of $V$. For $k=1$, the statement is just that of
Lemma~\ref{l.minimal-fully-unbounded}~. Hence, suppose that it holds for all
linear subspaces of dimension $< k$. Since the deviations are unbounded parallel
to all vectors $v\in V$, Lemma~\ref{l.minimal-fully-unbounded} and
Lemma~\ref{l.asymptotic-directions}(c) together imply that
\[
\sph^{d-1} \cap V \ \ssq \ \bigcup_{v\in {\cal A}^K_V} \overline{S^+(v)} \ .
\]
Hence ${\cal A}^K_V$ satisfies the assumptions of Lemma~\ref{l.hyperplanes}, and we
obtain vectors $v^1 \ld v^n$ and a linear subspace $\tilde V \ssq V$, such
that
\[
\sph^{d-1} \cap V \ \ssq \ \tilde V \cup \incup S^+(v^i) \ .
\]
By the induction assumption, there exists a residual set $\tilde{\cal R} \ssq
K$, such that for all $z\in\tilde{\cal R}$ the deviations parallel to all
$v\in\tilde V$ are unbounded both above and below. Further, the sets ${\cal
  R}_V(v^1) \ld {\cal R}_V(v^n)$ are residual by
Lemma~\ref{l.asymptotic-directions}(d). Hence ${\cal R} := \tilde{\cal R} \cap
{\cal R}_V(v_1) \cap \ldots \cap {\cal R}_V(v_n)$ is residual as the
intersection of a finite number of residual sets. For any $z\in{\cal R}$ and
$v\in \sph^{d-1}\cap V$ there holds $\sup_{n\in\N} D_v(n,z) = \infty$, since
either $v\in\tilde V$ or $v \in S^+(v^i)$ for some $i=1\ld n$. Obviously, this
also implies $\inf_{n\in\N} D_v(n,z) = - \sup_{n\in\N} D_{-v}(n,z) = - \infty$
for all $v\in V$. Thus ${\cal R}$ satisfies the assertions of the theorem.

\qed
\medskip

As the preceeding discussion shows, there is quite a good control about the
behaviour of the deviations $D(n,z)$ when these are unbounded and the dynamics
are minimal. In contrast to this, the Examples~\ref{e.unbounded} show that
similar results do not hold in general. The property which was exploited in the
construction of these examples was the existence of homotopically non-trivial
wandering domains, which is certainly a very strong form of
non-recurrence. Hence, an obvious question, which we leave open here, is whether
more can be said about systems with intermediate recurrence behaviour (apart
from Corollary~\ref{c.transitive-orbits}). For example, the following would be
interesting.
\begin{question}
  \label{q.unbounded-transitive} Suppose that $f \in\homd$ has lift $F$ with
  $\rho(F) \ssq \lambda v + \{v\}^\perp$ and the deviations parallel to $v$ are
  unbounded.  \alphlist
\item If $f$ is non-wandering (does not have any wandering open sets), does this
  imply that there is some $z\in\T^d$ with $\sup_{n\in\N} |D_v(n,z)| = \infty$?
\item If $f$ is topologically transitive, does this imply that there exists some
  $z\in\T^d$ with $\sup_{n\in\N} D_v(n,z) = \infty$, or even with $\sup_{n\in\N}
  D_v(n,z) = -\inf_{n\in\N} D_v(n,z) = \infty$?
\item Does the statement of Theorem~\ref{t.unbounded-deviations} remain true if
  $f$ is non-wandering/topologically transitive?
\item If $f$ is topologically transitive, or if $f$ is non-wandering and does
  not have periodic orbits, does $f$ have sensitive dependence on initial
  conditions? (In the non-wandering case, the existence of periodic orbits is
  excluded in order to avoid trivial counter-examples.)  \listend
\end{question}
A partial result on question (d) in dimension two will be given in the next
section. 
\medskip

Finally, as we have mentioned we will provide some examples with unbounded
deviations in all directions.
\begin{prop}
  \label{p.unbounded-examples}
  There exists a minimal pseudo-rotation $f\in\homd$, such that for any
  $v\in\R\smin\{0\}$ there holds $\sup_{n\in\N,z\in\T^d} |D_v(n,z)| = \infty$.
\end{prop}
\begin{bem} The proof below uses the well-known Anosov-Katok-method
  \cite{anosov/katok:1970,fayad/katok:2004} of fast approximation by rational
  maps. Using some standard modifications, one may additionally require $f$ to
  be uniquely ergodic and $\rho$ to be totally irrational. Further, one may
  require the examples to be smooth, but the rotation vectors will typically be
  Liouvillean. However, we expect that alternative methods should also yield
  similar examples with any given irrational rotation vector. The interesting
  question is whether such examples can also be made smooth, as this is the case
  for Herman's examples in \cite{herman:1983}.
\end{bem}

\begin{question}
  Are there smooth (${\cal C}^\infty$ or even real-analytic) irrational
  pseudo-rotations with Diophantine rotation vector and unbounded mean motion in
  all directions?
\end{question}
We will carry out the construction only in dimension $d=2$, again the
modifications for the general case are minor. Before we turn to the proof of
Proposition~\ref{p.unbounded-examples}, we need the following auxiliary
statement:
\begin{lem}
  \label{l.unbounded-criterion} Suppose $f\in \homeo_0(\torus)$ has lift $F$ and
  $\#\rho(F) = 1$. Further, assume that $u,w \in\R^d$ are linearly independent
  and for some $z\in\R^2$ there holds \romanlist
  \item $\sup_{n\in\N} |D_u(n,z)| = \sup_{n\in\Z} |D_w(n,z)| = \infty$;
  \item $\limsup_{n\to\infty} \frac{|D_u(n,z)|}{|D_w(n,z)|+1} = \infty$. 
  \listend
Then $\sup_{n\in\N}|D_v(n,z)| = \infty \ \forall v \in \R^2\smin 0$. 
\end{lem}
\proof Suppose for a contradiction that for some $v\in\R^2 \smin\{0\}$ and $C>0$
there holds $\sup_{n\in\N}|D_v(n,z)| \leq C$. Choose any vector $v' \in
\{v\}^\perp \smin\{0\}$, and write $u,w$ as $u=\lambda_1v+\lambda_2v'$ and
$w=\mu_1v+\mu_2v'$. In order to simplify notation, we will assume
$\lambda_1,\lambda_2,\mu_1,\mu_2\geq 0$. Note that due to property (i), there
holds $\lambda_2,\mu_2 \neq 0$. Now, given any $n\in\N$, we have
\begin{eqnarray*}
  \frac{|D_u(n,z)|}{|D_w(n,z)|+1} & = & \frac{|\lambda_1D_v(n,z) +
    \lambda_2D_{v'}(n,z)|}{|\mu_1D_v(n,z) + \mu_2D_{v'}(n,z)|+1} \\ 
  & \leq & \frac{\lambda_1 C + \lambda_2|D_{v'}(n,z)|}{\max\{1,\mu_2|D_{v'}(n,z)| - \mu_1C\}}
  \\ &\leq &  \max\left\{\lambda_1C+\frac{2\mu_1\lambda_2C}{\mu_2}, \frac{\lambda_1}{\mu_1} 
    + \frac{2\lambda_2}{\mu_2} \right\} \ .
\end{eqnarray*}
(For the last step, it is convenient to distinguished the cases $|D_{v'}(n,z)| <
\frac{2\mu_1C}{\mu_2}$ and $|D_{v'}(n,z)| \geq \frac{2\mu_1C}{\mu_2}$.) This
bound does not depend on $n$, which clearly contradicts (ii).

\qed
\medskip

\proof[Proof of Proposition~\ref{p.unbounded-examples}] For the construction, we
use the well-known Anosov-Katok-method. Before we start, we need some more
notation. Given $f,g\in \homeo_0(\torus)$, let $\dhom(f,g) :=
\sup_{z\in\torus}\{d(f(z),g(z)),d(f^{-1}(z),g^{-1}(z))\}$. If $F$ is a lift of
$f$, $n\in\Z$ and $z,v,\rho \in \R^d$, let $D_v(f,\rho,n,z) = \langle
F^n(z)-z-n\rho,v\rangle$ (hence, we do not suppress $f$ and $\rho$ in the
notation as usual, since we will consider sequences of maps and rotation vectors
below). Further, by $\lfloor x \rfloor$ we denote the integer part of
$x\in\R$.\medskip

By induction, we will choose sequences of rational rotation vectors $\alpha_n =
(p_n/q_n,p_n'/q_n)$, toral homeomorphisms $h_n$ and $f_n = h_n \circ
R_{\alpha_n} \circ h_n^{-1}$ and integers $k_n$ and $K_n$ with the following
properties:
\begin{list}{$(\roman{enumi})_n$}{\usecounter{enumi}}
\item $h_n(0) = 0$;
\item $|\alpha_n-\alpha_{n-1}| < 2^{-n}$;
\item $\dhom(f_n,f_{n-1}) < 2^{-n}$;
\item $\|(F^k_n(z)-z)/k-\alpha_j\| < 2^{-j} \quad \forall j=1\ld n,\ k=K_j \ld
  2K_j-1$;  
\item $|D_{e^1}(f_n,\alpha_n,k_j,0)| > j \quad \forall j=2,4\ld
  2\lfloor\frac{n}{2}\rfloor$;
\item $|D_{e^2}(f_n,\alpha_n,k_j,0)| > j\cdot (|D_{e^1}(f_n,\alpha_n,k_j,0)|+1)
  \quad \forall j=1,3\ld 2\lfloor\frac{n+1}{2}\rfloor -1$;
\item For all $z\in\torus$ and $j\leq n$, there holds $\bigcup_{k=0}^{q_j}
  B_{1/j}(f^k_n(z)) = \torus$.
\listend
%
%
Here $F_n$ denotes the lift of $f_n$. Before we carry out the induction, let us
see how $(i)$--$(vii)$ imply the statement of the proposition. Due to $(iii)$,
the limit $f=\nLim f_n \in\homeo_0(\torus)$ exists and is well-defined. By
$(ii)$, the limit $\alpha=\nLim \alpha_n$ exists, and there holds
$|\alpha-\alpha_n| < 2^{-n} \ \forall n\in\N$. From $(iv)$, we can therefore
deduce that for all $k \geq K_j$ and $z\in\T^d$ there holds
\[
\|(F^k(z)-z)/k-\alpha\| \ = \ \nLim \|(F^k_n(z)-z)/k-\alpha_j\| + 2^{-j} \ < \ 2^{-j+1}
\ .
\]
(Note that for all $n\geq j$ condition $(iv)_n$ yields $\|
(F^k_n(z)-z)/k-\alpha_j\| < 2^{-j} \ \forall k \geq K_j$.) This implies $\rho(F)
= \{\alpha\}$.

From $(v)$, it follows that $D_{e^1}(f,\rho,k_j,0) = \nLim
D_{e^1}(f_n,\alpha_n,k_j,0) \geq j \ \forall j \in 2\N$, in particular
$\sup_n|D_{e^1}(f,\rho,n,0)| = \infty$. Similarly, $\sup_n|D_{e^2}(f,\rho,n,0)|
= \infty$ follows from $(vi)$, which also implies
\[
\frac{|D_{e^2}(f,\rho,k_j,0)|}{|D_{e^1}(f,\rho,k_j,0)|+1} \ = \ \nLim
\frac{|D_{e^2}(f_n,\alpha_n,k_j,0)|}{|D_{e^1}(f_n,\alpha_n,k_j,0)|+1} \ \geq \ j
\quad \forall j \in 2\N+1 \ .
\]
Thus the assumptions of Lemma~\ref{l.unbounded-criterion} are all satisfied
(with $z=0$, $u=e^2$, and $w=e^1$), and we conclude that $f$ has unbounded
deviations in all directions. Finally, the fact that $f$ is minimal follows
easily from property $(vii)$.  \medskip

In order to start the induction, we choose $\alpha_0$ arbitrary and
$h_0=\Id_{\torus}$, such that $f_0=R_{\alpha_0}$. Then $(i)_0$ is trivial, and
$(ii)_0$--$(vii)_0$ are still void.

For the induction step $n\to n+1$, we certainly have to distinguish between odd
an even $n$, due to the asymmetric conditions $(v)_n$ and $(vi)_n$. However, we
will only give the details for even $n$, since the odd case is similar and even
easier. Thus, suppose $n$ is even and $\alpha_0\ld \alpha_n,\ h_0\ld h_n,\
k_0\ld k_n$ and $K_0 \ld K_n$ have been chosen and satisfy $(i)_j$--$(vii)_j \
\forall j=0\ld n$. Denote the lifts of $h_0\ld h_n$ by $H_0 \ld H_n$.

Let $C^j_n := \sup_{z\in\R^2,i=\pm 1} \pi_j(H^i_n(z)-z) \ (j=1,2)$, and note
that 
\[
\sup_{k\in\Z,z\in\R^2}|D_{e^1}(f_n,\alpha_n,k,z)| \leq 2C^1_l \ .
\]
Choose a $1/q_n$-periodic continuous function $\psi:\R\to\R$, which satisfies
$\psi(0)=0$ and
\[
   \psi(\theta) \ \geq \ (n+1)(2C^1_l+1)+C^2_n \quad \forall \theta \in I \ ,
\]
where $I \ssq (0,1/q_n)$ is an open and non-empty interval. Let $G(x,y) =
(x,y+\psi(x))$ and define $h_{n+1}$ via its lift $H_{n+1} = H_n\circ G$. Then
$(i)_{n+1}$ holds, since $G(0)=0$. Further, as $R_{\alpha_n} \circ G = G \circ
R_{\alpha_n}$, there holds $H_{n+1} \circ R_{\alpha_n} \circ H_{n+1}^{-1} = H_n
\circ R_{\alpha_n} \circ H_n = F_n$. Due to the uniform continuity of $H_{n+1}$
and its inverse, we can now ensure that conditions $(ii)_{n+1}$, $(iii)_{n+1}$
and $(v)_{n+1}$ hold just by choosing $\alpha_{n+1}$ sufficiently close to
$\alpha_n$. The same is true for all inequalities in $(iv)_{n+1}$ and
$(vi)_{n+1}$, except for the last ones with $j=n+1$, which do not yet make sense
since we have not chosen $k_{n+1}$ and $K_{n+1}$ so far. Similarly, we can
require that all but the last equalities in $(vii)_{n+1}$ hold (note that these
conditions are open with respect to $\dhom$.)

We will proceed in two steps and choose an intermediate rotation vector
$\tilde \alpha_{n+1}$ first. Suppose $q$ is large enough, such that
$\tilde\alpha_{n+1} = (p_n/q_n+1/qq_n,p_n'/q_n)$ is sufficiently close to
$\alpha_n$ in the above sense and $1/q < |I|$. Let $\tilde f_{n+1} = H_{n+1}
\circ R_{\tilde\alpha_{n+1}}\circ H_{n+1}^{-1}$.  Choose $k\in\N$, such that
$k/q \in I$ and let $k_{n+1} = kq_n$. We claim that
\begin{equation}
  \label{eq:7}
  |D_{e^2}(f_{n+1},\alpha_{n+1},k_{n+1},0)| \ > 
  \ (n+1) \ (|D_{e^1}(f_{n+1},\alpha_{n+1},k_{n+1},0)|+1) \ .
\end{equation}
In order to see this, note that the line $L= \kreis \times\{0\}$ is
$R_{\tilde\alpha_{n+1}}^{q_n}$-invariant, and
$R_{\tilde\alpha_{n+1}}^{k_{n+1}}(0) \in I \times\{0\}$ by the choice of $k$. Hence, we obtain that
\begin{eqnarray}
  \lefteqn{ D_{e^2}(\tilde f_{n+1},\tilde \alpha_{n+1},k_{n+1},0) \ = \ \pi_2\left(\tilde F_{n+1}^{k_{n+1}}(0)
      - k_{n+1}\tilde \alpha_{n+1}\right)} \nonumber \\ \label{eq:8}
  &\stackrel{(i)_{n+1}}{=}& \pi_2\left(H_{n+1}\circ R_{\tilde \alpha_{n+1}}^{k_{n+1}}(0)\right) - kp_n' 
  \ = \ \pi_2\left(H_n(k/q,\psi(k/q))+k(p_n,p_n')\right)- kp_n' \\
  &=& \pi_2\left(H_n(k/q,\psi(k/q)\right) \ \geq \ \psi(k/q)-C^2_n \ > \ (n+1)(2C^1_n+1) \ .
  \nonumber
\end{eqnarray}
On the other hand, note that $C^1_{n+1} = C^1_n$, since $G$ leaves the first
coordinate invariant. Therefore
\begin{equation}
  \label{eq:9}
  D_{e^1}(\tilde f_{n+1},\tilde \alpha_{n+1},k_{n+1},0) \ \leq \ 2C^1_N \ .
\end{equation}
Now (\ref{eq:8}) and (\ref{eq:9}) imply (\ref{eq:7}). 

In order to complete the induction step, note that due to uniform continuity
there exists $\delta > 0$, such that $d(z,z') < \delta$
implies $$d(H_{n+1}(z),H_{n+1}(z'))\ < \ \frac{1}{2(n+1)} \ . $$ We choose
$\alpha_{n+1}$ sufficiently close to $\tilde\alpha_{n+1}$, such that all
conditions in $(i)_{n+1}$, $(ii)_{n+1}$, $(iii)_{n+1}$, $(v)_{n+1}$ and
$(vi)_{n+1}$ are satisfied, as well as all except for the last ones in
$(iv)_{n+1}$ and $(vii)_{n+1}$. Further, we require that the orbits of
$R_{\alpha_{n+1}}$ are $\delta$-dense in \torus. Then the orbits of $f_{n+1} =
H_{n+1}\circ R_{\alpha_{n+1}}\circ H_{n+1}^{-1}$, which are periodic of period
$q_{n+1}$, will be $1/2(n+1)$-dense in \torus. This ensures the last
condition in $(vii)_{n+1}$. Finally, since $f_{n+1}$ is semi-conjugate to
$R_{\alpha_{n+1}}$ it has unique rotation vector $\alpha_{n+1}$, and hence for
any suitably large $K_{n+1}$ the last condition in $(iv)_{n+1}$ holds as
well. This completes the proof.

\qed

\section{Bounded mean motion and Lyapunov stability}

\label{LyapunovStability}

Suppose $X$ is a metric space and $f:X\to X$ is a continuous map. We call a
point $z\in X$ {\em $\eps$-Lyapunov stable}, if there exists $\delta>0$,
such that
\begin{equation}
  \label{eq:6}
f^n(B_\delta(z)) \ \ssq \  B_\eps(f^n(z)) \quad \forall n\in\N_0 \ .
\end{equation}
This is slightly weaker than the usual notion of a Lyapunov stable point, which
requires that for all $\eps>0$ there exists a $\delta>0$ with the property
(\ref{eq:6}).

$f$ has {\em sensitive dependence on initial conditions}, if there exists some
$\eps > 0$, such that for all $z\in X$ and $\delta>0$ there exist $z'\in
B_\delta(x)$ and $n\in\N_0$, such that $d(f^n(z),f^n(z'))) \geq \eps$. In this
situation, we call $\eps$ a {\em separation constant} for $f$.  It is easy to
see that $\eps$ is a separation constant for $f$ if and only if there does not
exist an $\eps$-Lyapunov stable point.

Now, assume that $f\in\homeo_0(\torus)$ has lift $F$, and let 
\begin{equation}
  \label{eq:7a}
  \eps_f \ := \ \sup\{\eps>0 \mid \forall z,z' \in \torus :  d(z,z')<\eps
   \follows d(f(z),f(z')) < 1/2\} \ .
\end{equation}
Note that for all $z,z' \in \R^2$ with $d(z,z') < \eps_f$ there holds
$d(F(z),F(z')) < \halb$, and therefore $d(F(z),F'(z)) = d(f(\pi(z)),f(\pi(z'))$.
In particular, a point $z\in\torus$ which is an $\eps_f$-Lyapunov stable point
for $f$ lifts to an $\eps_f$-Lyapunov stable point for $F$.

Given any point $z\in\R^2$, we call $\rho(F,z) = \nLim (F^n(z)-z)/n$ the
rotation vector of $z$, whenever the limit exists. We say $f$ is {\em
  non-wandering}, if there exists no open set $U\ssq\torus$ with $U\cap f^n(U) =
\emptyset \ \forall n\in\N$.
\begin{thm}
  \label{t.lyapunov-stability}
  Suppose $f\in\homeo_0(\torus)$ is non-wandering, and there exists an
  $\eps_f$-Lyapunov stable point $z_0$ with totally irrational rotation
  vector. Then $f$ is an irrational pseudo-rotation with bounded mean motion.
\end{thm}
As an immediate consequence, we obtain
\begin{cor} \label{c.sdic}
  Suppose $f\in\homeo_0(\torus)$ is a non-wandering irrational pseudo-rotation
  with unbounded mean motion. Then $f$ has sensitive dependence on initial
  conditions.
\end{cor}
We divide the proof of the theorem into the following two lemmas.
\begin{lem} \label{l.1}
   Suppose $f\in\homeo_0(\torus)$ is non-wandering, and there exists an
  $\eps_f$-Lyapunov stable point $z_0$ with totally irrational rotation
  vector. Then $f$ is an irrational pseudo-rotation. 
\end{lem}
\proof Let $F$ be a lift of $f$ and suppose $z_0$ is an $\eps_f$-Lyapunov stable
point, $\hat z_0\in\pi^{-1}(z)$ and
\begin{equation} \label{eq:11} \nLim (F^n(\hat z_0)-\hat z_0)/n \ = \ \rho \ \notin \  \Q^2 \ .
\end{equation}
As mentioned above, $\hat z_0$ is an $\eps_f$-Lyapunov stable point for $F$. We
assume without loss of generality that $\hat z_0=0$. Choose $\delta > 0$, such
that $F^n(B_\delta(0)) \ssq B_{\eps_f}(F^n(0)) \ \forall n \in \N_0$. Let
$U=B_\delta(0)$ and $U_0=\pi(U)$. Since $\rho$ is not rational, $U$ does not
contain any periodic points. The fact that $f$ is non-wandering therefore
implies that there exist infinitely many pairs $(m,v) \in \ \N\times \Z^2$, such
that
\begin{equation} \label{eq:13}
F^m(U) \cap (U+v) \ \neq \ \emptyset \ .
\end{equation}
Given any such pair, let $\tilde F := F^m-v$ and $\tilde \rho = m\rho -v$, such
that
\begin{equation} \label{eq:13a} \nLim \tilde F^n(0)/n \ = \ \tilde \rho \ .
\end{equation}

Without loss of generality, we assume $\pi_1(\tilde\rho) > 0$. Our aim is to
show that $\rho(\tilde F) \ssq \tilde\rho\R$, this will then imply the statement
of the lemma quite easily. In order to do so, let $V := \bigcup_{n\in\N_0}
\tilde F^n(U)$. As $\tilde F(U) \cap U \neq \emptyset$, the set $V$ is
connected. Further, (\ref{eq:13a}) together with the Lyapunov stability of $\hat
z=0$ implies that for each fixed $k\in\N_0$ there holds
\begin{equation}
  \label{eq:17} 
  \nLim \ntel d\left(\tilde F^n(U)+(0,k),n\tilde\rho\right) \ = \ 0 \ .
\end{equation}
We claim that for sufficiently large $N\in\N$, the integer translate $V+(0,N)$
is disjoint from $V$. In order to see this, note that due to (\ref{eq:17}) the
sets $U$ and $V+(0,N)$ will be disjoint for sufficiently large $N$, since only a
finite number of iterates of $U+(0,N)$ intersect the strip $[-1,1]\times \R
\supseteq U$. Hence, if the orbit of $U$ intersects $V+(0,N)$, then it must
first intersect $U+(0,N)$. However, by the same argument the orbit of $U$ can
only intersect a finite number of its vertical integer translates, such that for
sufficiently large $N$ we have $V\cap V+(0,N) = \emptyset$ as required.

Now let $a=\inf\{x\in\R \mid (0,x) \in V\}$, $b=\sup\{x\in\R \mid (0,x) \in
V+(0,N)\}$ and define $W$ as the union of $V$, $V+(0,N)$ and the vertical 
arc from $(0,a)$ to $(0,b)$. Let $Y$ be the unique connected component of
$\R^2\smin \overline{W}$ which is unbounded to the left, and $A = \R^2 \smin
Y$. (See Figure~\ref{f.1}.) The following two remarks about these objects will be
helpful below:

First, from (\ref{eq:17}) we can deduce for any $m_1<\rho_2/\rho_1 < m_2$ there
exist constants $c_1,c_2$, such that
\begin{equation}
   \label{eq:18}
   A \ \ssq \ \{ (x,y) \in \R \mid m_1x+c_1 \leq y \leq m_2x+c_2 \} \ .
\end{equation}
Secondly, due to the definition of $W$, its connectedness and the fact that it
`stretches out to infinity' by (\ref{eq:17}), the set $S := (\R^+ \times \R)
\smin A = (\R^+ \times \R) \cap Y$ consists of exactly two connected
components. These can be defined as follows: Fix any $\zeta_0 \in Y$ with
$\pi_1(\zeta) < 0$. For any $\zeta \in S$, there is a path $\gamma_\zeta$ from
$\zeta$ to $\zeta_0$. Let $x_\zeta$ be the second coordinate of the first point
in which $\gamma_\zeta$ intersects the vertical axis. The fact whether $x_\zeta$
lies above $b$ or below $a$ does not depend on the choice of the path, since
this would contradict the connectedness of $W$. It is now easy to see that $S^-
= \{\zeta \in S \mid x_\zeta > b\}$ and $S^+ = \{ \zeta \in S \mid x_\zeta <
a\}$ form a partition of $S$ into two connected components.

\begin{figure}
\begin{minipage}[t]{\linewidth}  
  \epsfig{file=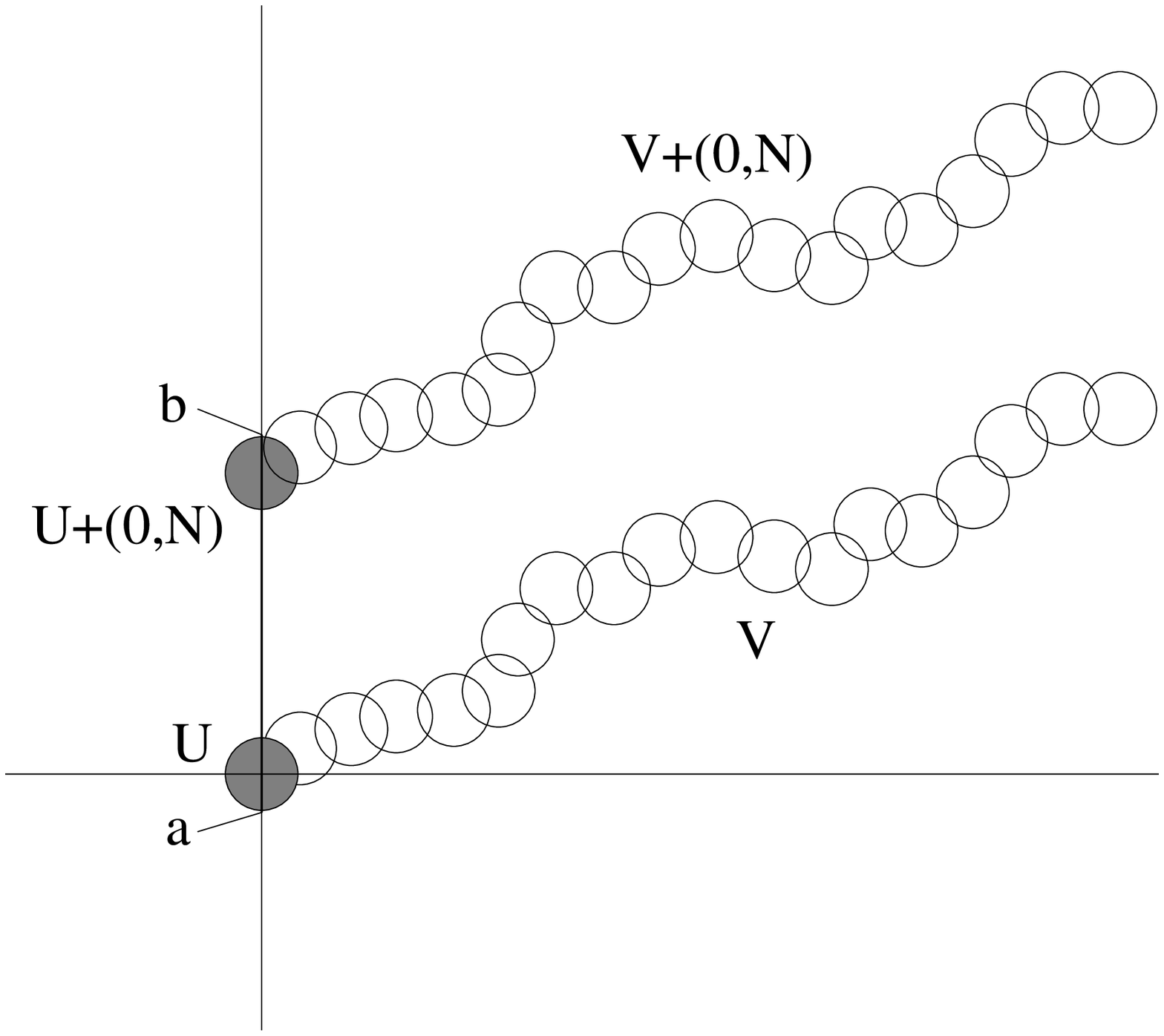, clip=, width=0.6\linewidth} \hspace{-6eM}
  \epsfig{file=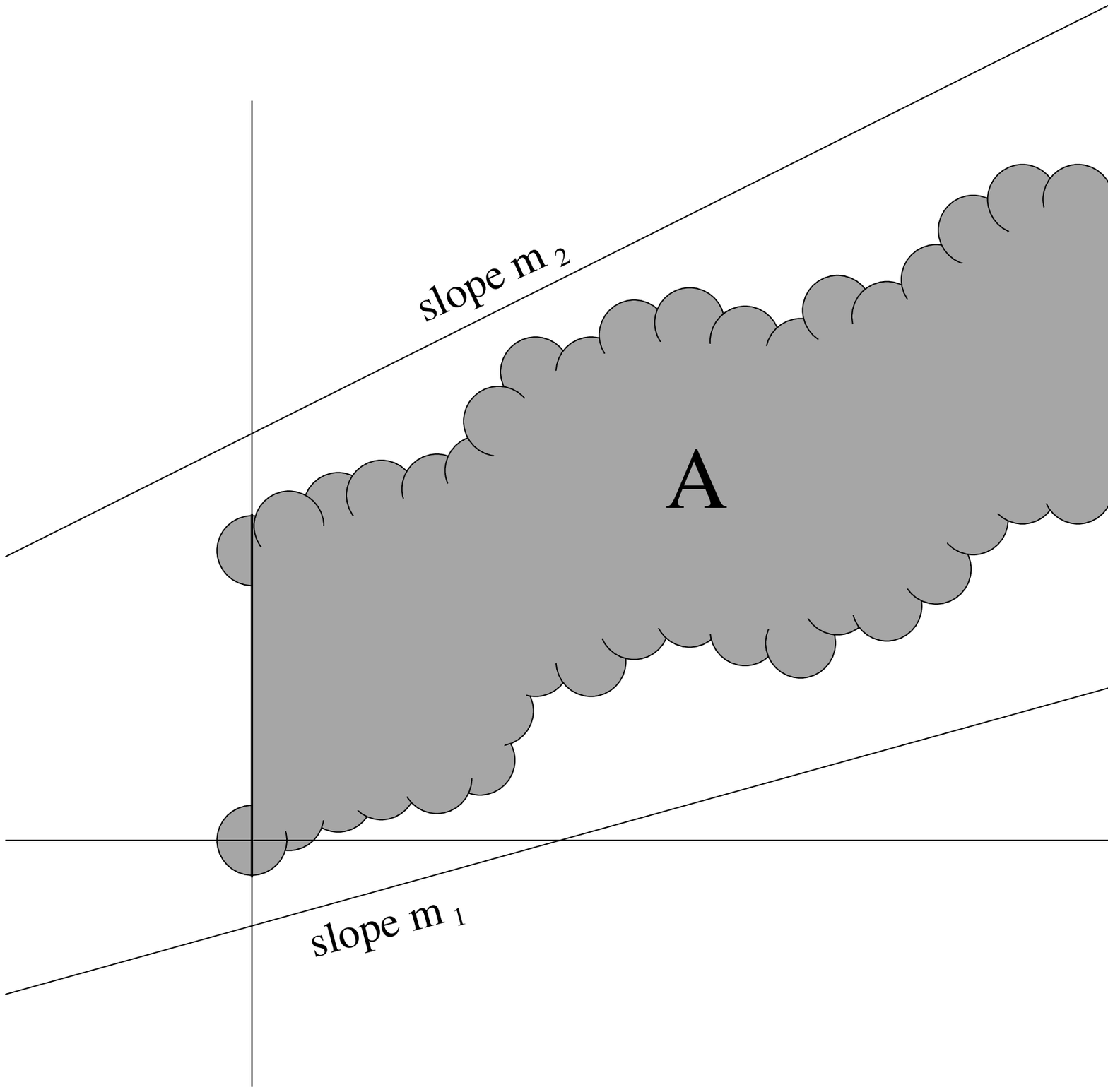, clip=, width=0.6\linewidth}
  \caption{\small The construction of the sets $V$ and $W$ on the right, and the
    set $A$ on the left.  \label{f.1} }
        
\end{minipage}
\end{figure} 

\begin{claim} \label{c.rho} There exists some constant $\alpha>0$, such that
  $z\in A$ and $\pi_1(z) \geq \alpha$ imply $\tilde F(z) \in A$ and $\tilde F^{-1}(z) \in
  A$
\end{claim}
\proof We show that there exists a constant $\alpha' > 0$, such that $\pi_1(z)
\geq \alpha'$ and $z\in S$ implies $\tilde F(z) \in S$ and $\tilde F^{-1}(z)\in
S$. If we let
\begin{equation} \label{e.M}
M\ :=\ \sup_{z\in\R^2}\|\tilde F(z)-z\| \ = \ \sup_{n\in\R^2} \|\tilde F^{-1}(z)-z\| \ ,
\end{equation}
then for any $z\in\pi_1^{-1}[\alpha'+M,\infty)$ this implies that there holds
$\tilde F(z) \in S \follows z\in S$, and hence $z\in A \follows \tilde F(z) \in
A$. The same applies to $\tilde F^{-1}(z)$. Thus we can choose
$\alpha=\alpha'+M$.

From (\ref{eq:17}) we deduce that there exists some $n_0\in\N$, such that
$$\pi_1\circ \tilde F^n(U) \ \ssq \ (4M,\infty) \quad \forall n\geq n_0 \ . $$ Choose
$\alpha'>0$, such that for all $j\leq n_0$ there holds $\pi_1\circ \tilde F^j(U)
\ssq [0,\alpha')$. We obtain the following statement:
\begin{equation}
  \label{eq:19}
  \pi_1\circ\tilde F^n(U) \cap [\alpha',\infty) \neq \emptyset \quad
  \follows \quad \pi_1\circ\tilde F^k(U) \cap [0,4M] = \emptyset \quad  \forall k\geq n \ .
\end{equation}
Of course, the same statement applies to $U+(0,N)$. 

Due to (\ref{eq:18}), it is possible to choose some $K>0$, such that 
\[
B \ := \ [0,4M] \times [K-M,\infty) \ssq \ S^+ \ .
\]
Let $z^*=(3M,K)$ and fix $z \in S^+$ with $\pi_1(z) \geq \alpha'$. Then, since
$S^+$ is open and connected, there exists a simple path $\gamma : [0,1] \to S^+$
from $z$ to $z^*$.  We claim that $\gamma$ can be chosen such that its image
contained in $S^+ \cap \pi^{-1}[3M,\infty)$.

Suppose not, and let $t_0 := \min\{t\in[0,1] \mid \gamma(t) \in B\}$ and $\Gamma
= \{\gamma(t) \mid t\in[0,t_0]\}$. Then $\Gamma$ divides the set $([0,\pi_1(z)]
\times \R) \smin B$ into exactly two connected components $D^+$ and $D^-$, which
are unbounded above, respectively below (see Figure~\ref{f.2}). Now, if
$\overline{W}$ does not intersect $D^+ \cap \pi^{-1}_1[0,4M)$, then $D^+\cap
\pi^{-1}_1[0,4M) \ssq S^+$, and it is easy to see that in this case either
$\gamma$ does not intersect $\pi_1^{-1}[0,3M)$, or we can modify it to that
end. Otherwise, there must be some $k\in\N$, such that $F^k(U)$ or
$F^k(U)+(0,N)$ intersects $D^+\cap \pi^{-1}_1[0,4M)$. However, since $U$
intersects $D^-$ and the set $\bigcup_{i=0}^k \tilde F^i(U)$ is connected, this
implies that there must be some $n\leq k$, such that $\pi_1\circ\tilde F^n(U)$
intersects $[\pi_1(z),\infty] \ssq [\alpha',\infty]$. This contradicts
(\ref{eq:19}).

Summarising, we have found a path $\gamma$ from $z$ to $z^*$, which is contained
in $S^+ \cap \pi^{-1}_1[3M,\infty)$. In particular, $\gamma$ is contained in the
complement of $\overline W$. Consequently, the path $\tilde F\circ \gamma$ is
contained in the complement of $\tilde F(\overline{W})$. At the same time, it is
also contained in $\pi^{-1}_1[2M,\infty)$. However, it follows from the
construction of $W$ and the definition of $M$ that $\tilde F(\overline W) \cap
\pi^{-1}_1[2M,\infty) = \overline W \cap \pi^{-1}_1[2M,\infty)$. Hence, the path
$\tilde F \circ \gamma$ is also contained in the complement of $\overline{W}$.
Furthermore, it joins $\tilde F(z)$ to the point $\tilde F(z^*)$, which is
contained in $B \ssq S^+$. This implies that $\tilde F(z)$ is equally contained
in $S^+$. When $z\in S^-$, the argument is similar. In the same way, one can
show that $\tilde F^{-1}(z) \in S$, which completes the proof.

\roundqed

\begin{figure} 
\begin{minipage}[h]{\linewidth}  
\begin{center}
  \epsfig{file=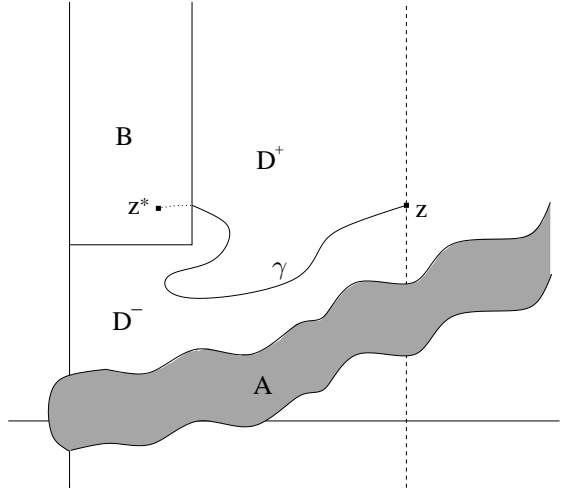, clip=, width=0.5\linewidth}
\caption{\small The path $\gamma$ divides $([0,\pi_1(z)]\times \R) \smin B$
    into two connected components $D^+$ and $D^+$. As argued in the proof, the
    set $W$ (and hence $A$) cannot `reach back' into $D^+$ too far.  \label{f.2}}
 \vspace{-29ex} $\gamma$ \vspace{27ex}
 \end{center}
\end{minipage}
\end{figure} 

Using this `almost-invariance' of $A$, we can now prove the following.
\begin{claim} \label{cl.rotationset}
  There holds $\rho(\tilde F) \ssq \tilde\rho\R$.
\end{claim}
\proof Suppose for a contradiction that $\sigma \in \rho(\tilde F) \smin
\tilde\rho\R$. Without loss of generality, we may assume that $\sigma$ is an
extremal point of $\rho(\tilde F)$. In this case, $\sigma$ is realised by some
orbit \cite{misiurewicz/ziemian:1989}, that is, there exists some $z_1 \in \R^2$
with $\lim_{|n|\to \infty} (\tilde F^n(z_1)-z_1)/n = \sigma$. Note that the
limit is take two-sided.

There holds $\sigma \neq 0$. Hence, by permuting the coordinates if necessary,
we may assume $\pi_1(\sigma) \neq 0$. We first assume that $\pi_1(\sigma) >
0$. In this case, there exists some $r\in\R$, such that $\pi_1(\tilde
F^n(z_1)-z_1) > -r \ \forall n\in\N_0$. Now, there holds $\R^+ \times \R \ssq
\bigcup_{j\in\N} A+(0,jN)$. Replacing $z_1$ by an integer translate if
necessary, we may therefore assume $z_1\in A$ and $\pi_1(z_1) \geq
\alpha+r$. In this case, Claim~\ref{c.rho} implies $\tilde F^n(z_1) \in A \
\forall n\in\N_0$, and therefore (\ref{eq:18}) yields $\sigma \in \tilde\rho\R$,
contradicting our assumption.  When $\pi_1(\sigma)<0$, we can proceed in the
same way by considering the backwards orbit of $z_1$.


\roundqed
\medskip

Now we can conclude the proof of the lemma. Due to (\ref{eq:11}), and using the
fact that $\rho$ is totally irrational and $v$ is an integer vector, there
exists some $n_0 \in \N$, such that for all $n\geq n_0$ there holds
$d(F^n(z)-z,\R v) > 2\eps_f$. It follows that $F^n(U) \cap (U+\lambda v) =
\emptyset \ \forall n \geq n_0,\ \lambda \in \R$. Hence, we can find a second
pair $(m',v')$ that satisfies (\ref{eq:13}), such that the vectors $v$ and $v'$
are linearly independent. Now, applied to $(m,v)$, Claim~\ref{c.rho} yields
\begin{equation}
  \label{eq:14}
  \rho(F) \ = \ \frac{\rho(\tilde F)}{m} + \frac{v}{m} \ \ssq \ 
  \left(\rho-\frac{v}{m}\right) \R + \frac{v}{m} \ . 
\end{equation}
Similarly, we obtain 
\begin{equation}
  \label{eq:15}
  \rho(F)  \ \ssq \  \left(\rho-\frac{v'}{m'}\right) \R + \frac{v'}{m'}
\end{equation}
for the pair $(m',v')$. As $\rho$ is totally irrational and the rational vectors
$v/m$ and $v'/m'$ are linearly independent, the vectors $\rho-v/m$ and
$\rho-v'/m'$ are linearly independent as well. Thus (\ref{eq:14}) and
(\ref{eq:15}) together imply that $\rho(F)$ is reduced to a singleton. Since the
rotation vector $\rho$ is realised on the orbit of $z_0$, we must have $\rho(F)
= \{\rho\}$.

\qed

\begin{lem} \label{l.2}
  Suppose $f\in\homeo_0(\torus)$ is a non-wandering irrational pseudo-rotation,
  and there exists an $\eps_f$-Lyapunov stable point $z_0$. Then $f$ has bounded
  mean motion.
\end{lem}
\proof We use the same definitions and notation as in the proof of the
preceeding lemma. As in the proof of Proposition~\ref{p.unbounded-examples}, we
write $D_v(f,\rho,n,z)$ for the deviations from the constant rotation, in order
to emphasize the dependence both on the map and on the rotation vector. First,
let $u$ be orthogonal to $\rho-v/m$ with norm $\|u\|=1$, where the pair $(v,m)$
is chosen as in (\ref{eq:13}). Suppose that the deviations $D_u(f,\rho,n,z_0)$
are not uniformly bounded in $n\in\N_0$, without loss of generality
\begin{equation} \label{eq:20}
  \sup_{n\in\N_0} D_u(f,\rho,n,z_0) \ = \ \infty \ .
\end{equation}
In this case, an easy computation yields that the deviations $D_u(\tilde f,\tilde
\rho,n,z_0)$ are equally unbounded in $n\in\N_0$. (By $\tilde f$ we denote the
toral homeomorphism induced by $\tilde F$.) We want to lead this to a
contradiction.

Assume without loss of generality that $\pi_1(\tilde \rho) > 0$. Note that since
$\tilde f$ is an irrational pseudo-rotation, the quantities $(\tilde
F^n(z)-z)/n$ converge to $\tilde \rho$ uniformly on $\R^2$. It follows that
there exists some $n_0\in\N$, such that there holds
\begin{equation} \label{eq:12}
  \pi_1(\tilde F^n(z)-z) \ \geq \ M+1 \quad \forall n\geq n_0,\ z\in\R^2 \  ,
\end{equation}
where $M$ is defined as in (\ref{e.M}).

Together with the Lyapunov stability of $z_0$ and its lift, the fact that the
deviations $D_u(\tilde f,\tilde \rho,n,z_0)$ are unbounded implies that for
every $r>0$ there exists some $n\in\N_0$, such that $$\langle z, u\rangle \ \geq
\ r \quad \forall z \in \tilde F^n(U) \ . $$ It follows that for all $k \in
\{n-n_0\ld n+n_0\}$ and all $z \in \tilde F^k(U)$ there holds $\langle
z,u\rangle \geq r-n_0 M$. By (\ref{eq:12}), we obtain that there exists an
interval $I$ of length $2M$, such that all $z \in A \cap \pi^{-1}_1(I)$ satisfy
$\langle z,u\rangle \geq r-n_0 M - N$. Furthermore, we may assume that $I$ lies
to the right of the point $n_0 M+1$.

Now, as argued in the proof of Lemma~\ref{l.1}, any point in $\T^2$ has a lift
in the compact set $A_0 := A \cap \pi^{-1}_1(n_0 M,n_0 M+1]$. Furthermore, if
$n_0$ is chosen sufficiently large, such that $n_0M \geq \alpha$ with $\alpha$
as in Claim~\ref{c.rho}, then the forward orbit of this lift under $\tilde F$
will stay in $A$ all the time. However, any such orbit has to pass though the
set $\pi^{-1}_1(I)$. By choosing $r \geq n_0 M + N + \sup_{z\in A_0}\|z\| +1$,
we obtain that for every $z\in\torus$ there holds $\sup_{n\in\N} D_u(\tilde
f,\tilde\rho,n,z) \geq 1$. This contradicts Lemma~\ref{l.semi-bounded-orbits}~.

Hence, we may assume that $z_0$ has bounded deviations parallel to $u$. It
follows that that the set $A$ lies within a bounded distance of the semi-line
$\tilde \rho \R$. Since every point in $\torus$ has a lift whose orbit stays
forever in $A$, this further implies that the deviations $D_u(\tilde f,\tilde
\rho,n,z)$ are uniformly bounded in $n$ and $z$. Consequently, the same holds
for the deviations $D_u(f,\rho,n,z)$.

Now replace the pair $(m,v)$ by $(m',v')$ with $v$ and $v'$ linearly
independent, as at the end of the proof of Lemma. Then we obtain another,
linearly independent vector $u'$, such that the deviations $D_{u'}(f,\rho,n,z)$
are uniformly bounded. This implies that the vectors $D(f,\rho,n,z)$ are
uniformly bounded as well.

\qed

\begin{bem}
  The only time that we used the fact that the rotation vector $\rho$ of $z_0$
  is totally irrational, and not only non-rational, in the above proof was to
  deduce that there exist two pairs $(m,v)$ and $(m',v')$ satisfying
  (\ref{eq:13}), with linearly independent vectors $v$ and $v'$. Hence, if
  $\rho$ is neither rational nor totally irrational, the above construction
  still yields the following statement: \smallskip

  Suppose $f\in\homeo_0(\torus)$ is non-wandering, and there exists an
  $\eps_f$-Lyapunov stable point $z_0$ with a rotation vector $\rho$ that is
  neither rational nor totally irrational. Then there exists a vector
  $w\in\Q^2$, such that $\rho(f) \ssq \rho + (\rho-w)\R$ and the deviations
  orthogonal to $\rho-w$ are uniformly bounded.
\end{bem}



\bibliography{qpfs,torus,dynamics} \bibliographystyle{unsrt}

\end{document}